\documentclass[12pt]{article}
\usepackage{amssymb}
\usepackage{amsmath}

\newcommand\F{{\cal F}}
\newcommand\R{{\cal R}}

\newcommand\G{{\cal G}}
\newcommand\C{{\cal C}}

\newcommand\qed{\hbox{}\hfill $\Box$}

\newtheorem{thm}{Theorem}
\newtheorem{cor}{Corollary}
\newtheorem{conj}{Conjecture}

\newtheorem{lem}{Lemma}
\newtheorem{rem}{Remark}

\title{A general 2-part Erd\H os-Ko-Rado theorem}

\author{
{\bf Gyula O.H. Katona}\thanks{This research was supported
by the National Research, Development and
Innovation Office -- NKFIH Fund No's SSN117879 and K116769.}\\
MTA R\'enyi Institute\\
Budapest Pf 127, 1364 Hungary \\ ohkatona@renyi.hu }

\begin{document}
\date{}

\maketitle

\begin{abstract}
A two-part extension of the famous Erd\H{o}s-Ko-Rado Theorem is proved. The underlying set is partitioned into
$X_1$ and $X_2$. Some positive integers $k_i, \ell_i (1\leq i\leq m)$
are given.  We prove that if $\F$ is an intersecting family containing members $F$ such that $|F\cap X_1|=k_i, |F\cap X_2|=\ell_i$ holds for one of the values $i (1\leq i\leq m)$ then $|\F|$ cannot exceed the size of the largest subfamily containing one element.

\end{abstract}

\section{Introduction}
Let $X$ be a finite set of $n$ elements. A family $\F\subset 2^{X}$
is called {\it intersecting} if $F,G\in \F$ implies $F\cap G\not= \emptyset$. The family of all $k$-element subsets of $X$ is denoted by ${X\choose k}$. The celebrated theorem of Erd\H os, Ko and Rado  is the following.
\begin{thm} {\rm{\cite{EKR}}} Suppose that an integer $k\leq {n\over 2}$ is given and $\F\subset {X\choose k}$ is intersecting. Then
$$|\F|\leq {n-1\choose k-1}.$$
\end{thm}
The family of all $k$-element subsets containing a fixed element $x\in X$ shows that the estimate is sharp.

The goal of our paper is to consider the problem when the underlying set is partitioned into two parts $X_1, X_2$ and the sets $F\in \F$
have fixed sizes in both parts. More precisely
let $X_1$ and $X_2$ be disjoint sets of $n_1$, respectively $n_2$ elements. \cite{F} considered
such subsets of $X=X_1\cup X_2$ which had $k$ elements in $X_1$ and $\ell$ elements in $X_2$. The family of all such sets is denoted by ${X_1, X_2\choose k, \ell}$.
The construction above, taking all possible sets containing a fixed element also works here. If the fixed element is in $X_1$
then the number of these sets is
$$ {n_1-1\choose k-1}{n_2\choose \ell},$$
otherwise it is
$${n_1\choose k}{n_2-1\choose \ell-1}.$$
The following theorem of Frankl \cite{F} claims that the larger one of these is the best.
\begin{thm} Let $X_1, X_2$ be two disjoint sets of $n_1$ and $n_2$ elements, respectively. The positive integers $k, \ell$ satisfy the inequalities $2k\leq n_1, 2\ell \leq n_2$. If $\F$ is an intersecting subfamily of ${X_1, X_2\choose k, \ell}$ then
$$|\F|\leq \max \left\{  {n_1-1\choose k-1}{n_2\choose \ell},
{n_1\choose k}{n_2-1\choose \ell-1}\right\}.$$
\end{thm}

The goal of the present paper is to generalize Theorem 2 for the case when other sizes are also allowed that is the family consists of sets
satisfying $|F\cap X_1|=k_i, |F\cap X_2|=\ell_i$ for certain pairs of integers. Using the notation above, we will
consider subfamilies of
$$\bigcup_{i=1}^m{X_1,X_2\choose k_i,\ell_i}.$$

The generalization is however a little weaker at one point. In Theorem 2 the thresholds $2k\leq n_1, 2\ell \leq n_2$ for validity are natural. If either $n_1$ or $n_2$ is smaller then the problem becomes trivial, all such sets can be selected in $\F$. In the generalization below there is no such natural threshold. There will be another difference in the formulation. We give the construction of the extremal family rather than the maximum number of sets. A family is called {\it trivially intersecting} if there is an element contained in every member.

\begin{thm} Let $X_1, X_2$ be two disjoint sets of $n_1$ and $n_2$ elements, respectively. Some positive integers $k_i, \ell_i (1\leq i\leq m)$
are given. Define $b=\max_i \{ k_i, \ell_i\}$. Suppose that $9b^2\leq n_1, n_2$. If $\F$ is an intersecting subfamily of
$$\bigcup_{i=1}^m{X_1,X_2\choose k_i,\ell_i}$$
then $|\F|$ cannot exceed the size of the largest trivially intersecting family satisfying the conditions.
\end{thm}

Section 2 gives the proof of Theorem 3, using the method of ``cyclic permutations". Section 3 contains some open questions, while Section 4 shows some similar results from the past as motivations.

\section{The proof of the main theorem}

We will use the method of cyclic permutations \cite{Ksp} giving a simple
proof of the EKR theorem. There the analogous problem is solved for
intervals along a cyclic permutation and then a double counting easily
finishes the proof. Here we need a pair of cyclic permutations: one for
$X_1$ and one for $X_2$. A cycle of size $n_i$ will be represented by the
integers $\mod n_i$. The usual notation is $\mathbb{Z}_{n_i}$. Hence the
pair of cycles will be $\mathbb{Z}_{n_1}\times \mathbb{Z}_{n_2}$. The
direct product of the intervals of length $k$ and $\ell$, in $\mathbb{Z}_{n_1}$ and $\mathbb{Z}_{n_2}$, respectively,
will be a $k\times \ell$ {\it rectangle}
in $\mathbb{Z}_{n_1}\times \mathbb{Z}_{n_2}$. Problems analogous to our Theorem 3 will be considered
for such rectangles.

The proof will be divided into lemmas. If only one cycle is involved we
use the notation $\mathbb{Z}_n$ or even forget about this notation.

The {\it distance} $d(u,v)$ of the elements $u,v\in \mathbb{Z}_n$ is
the smaller distance along the cycle.

\begin{lem} {\rm{(Folklore)}}  Suppose $2\leq 2k< n$ and let  $G_n^k=(\mathbb{Z}_n ,E)$ be a graph where two vertices are adjacent if their distance (mod $n$) is at most $k-1$. Then the largest click (complete subgraph) in $G_n^k$ has $k$ vertices. If $k$ vertices form a click then
these vertices are consecutive (mod $n$).
\end{lem}

{\bf Proof.} $k$ consecutive elements obviously form a click. In order to prove that there is no larger click, choose a vertex of a click. By symmetry
one can suppose that $0$ is this vertex. The potential other vertices of the click are $-(k-1), -(k-2), \ldots , -1, 1, \ldots , k-2, k-1$.
By the assumption $2k< n$ all these vertices are distinct. The pairs
$(-(k-1), 1), (-(k-2), 2), \ldots , (-1, k-1)$ are not adjacent, therefore the click can contain only one of them, proving that the click
contains at most $1+(k-1)=k$ vertices.

To prove the second statement of the lemma suppose that $U\subset \mathbb{Z}_n$ spans a click, moreover $|U|=k, 0\in U$ hold. Let $i\in U (1<i\leq k-1)$ be in the click. Its ``pair" cannot be chosen: $-(k-i)\not\in U$. However the distance of $i$ and $-(k-i-1)$ is also at least $n-(k+1)\geq 2k+1-(k+1)=k$. Therefore  $-(k-i-1)\not\in U$ also holds. Its ``pair" $i-1$ must be in $U$. Let $j$ be the largest element of $U$. We have seen that $j-1, j-2, \ldots , 1$ are in $U$. On the other hand, since $j+1, \ldots , k-1\not\in U$, their ``pairs"
must be in $U$. It is fully determined: $U=\{ -(k-j-1), -(k-j-2), \ldots , 0, 1, \ldots , j\}$.
 \qed

\smallskip
An {\it interval of length} $a$ in $\mathbb{Z}$ is a set of form $\{ i+1, i+2, \ldots , i+a\}$ (mod $n$). The {\it left end} of this interval is $i+1$. The distance of two intervals is the minimum distance between two elements, one from each intervals.

\begin{lem} Let $k, b, n$ be positive integers satisfying $2(k+b)\leq n$. If $k+b+1$ distinct intervals of length $k$ are given in $\mathbb{Z}_n$ then there is a pair among them whose distance is at least $b+1$.
\end{lem}

{\bf Proof.} Apply Lemma 1 with  $k+b$. The left ends of our intervals cannot form a click in $G_n^{k+b}$ therefore two of them must have a distance at least $k+b$. Consequently the distance of these intervals is at least $b+1$. \qed

\smallskip
Let $I$ be an interval of length $k$ in $\mathbb{Z}_{n_1}$ while $J$ is an interval of length $\ell$ in $\mathbb{Z}_{n_2}$. The direct product $I\times J$ is a $k\times \ell$ {\it rectangle} in $\mathbb{Z}_{n_1}\times \mathbb{Z}_{n_2}$.
We say that the rectangles $I_1\times J_1$ and $I_2\times J_2$ are {\it proj-intersecting} if either $I_1\cap I_2$ or $J_1\cap J_2$ is non-empty. A family of rectangles is {\it proj-intersecting} if any two rectangles in the family are proj-intersecting. If $\R$ is a family of rectangles,
let $\R^1$ denote the set of intervals obtained by projecting the members of $\R$ on $\mathbb{Z}_{n_1}$. The family $\R^2$ is defined similarly. The inequality
$$|\R|\leq |\R^1|\cdot |\R^2|\eqno(1)$$
is obvious.

\begin{lem} Suppose that the positive integers $k, \ell, b, n_1, n_2$ satisfy the inequalities $k, \ell \leq b, 2(k+b)\leq n_1, 2(\ell+b)\leq n_2$.
Let $\R$ be a proj-intersecting family of size $|\R|\geq 9b^2$ of $k\times \ell$  rectangles  in $\mathbb{Z}_{n_1}\times \mathbb{Z}_{n_2}$. Then
either there are two rectangles $R_1, R_2\in \R$ such that
$$R_1=I_1\times J_0, R_2=I_2\times J_0, d(I_1,I_2)\geq b+1 \eqno(2)$$
or there are two rectangles $R_3, R_4\in \R$ such that
$$R_3=I_0\times J_3, R_4=I_0\times J_4, d(J_3, J_4)\geq b+1.\eqno(3)$$
\end{lem}

{\bf Proof.} (1) implies that either $|\R^1|>2b$ or $|\R^2|>2b .$

First suppose that $|\R^1|>2b\geq k+b$. Using Lemma 2 two members of $\R^1$ are obtained such that the distance between these two intervals $I_5, I_6$  is at least $b$. Let $R_5$ and $R_6$ be two rectangles whose projections to $\mathbb{Z}_{n_1}$
are $I_5, I_6$, respectively. Of course $I_5\cap I_6$ is empty, moreover no interval of length $k$ can intersect both of them, since $k\leq b$. Choose an arbitrary $R_7\in \R$ with projection $I_7$. It can intersect at most one of $I_5$ and $I_6$.  Since $\R$ is proj-intersecting, $J_5\cap J_6\not= \emptyset$ and one of $J_7\cap J_5, J_7\cap J_6$ is non-empty, too. The length of $J_5\cup J_6$ is at most $2\ell-1$. Here $J_7$ must meet $J_5\cup J_6$ therefore there are at most $3\ell -2$ choices for
$J_7$. Hence we have
$$|\R^2|\leq 3\ell -2< 3b.\eqno(4)$$

The {\it multiplicity} $\mu (J)$ of $J\in \R^2$ is the number of members of $\R$ with projection $J$ on $\mathbb{Z}_{n_2}$. Obviously $\sum_{J\in \R^2}\mu (J)=|\R|$ holds.  $|\R|\geq 9b^2$ and (4) imply that there is a member $J_0$ of $\R^2$ having multiplicity at least $3b\geq k+b+1$.  Lemma 2 implies again that there are two rectangles $R_1, R_2\in \R$ such that they both have the same projection $J_0$ on $\mathbb{Z}_{n_2}$ while their projections $I_1, I_2$ have distance at least $b+1$. A pair of rectangles of form (2) was found.

The other case when $|\R^2|>2b $ is analogous, then a pair of type (3) can be found. \qed

\smallskip
We say that the rectangles $R_1$ and $R_2$ in (2) form a $b${\it -blocking pair} with {\it base} $J_0$.

\smallskip
\begin{lem} Suppose that the rectangles $R_1, R_2$ form a $b$-blocking pair with base $J_0$ and these two and a third rectangle $U\times V$ are pairwise proj-intersecting where
$U\subset \mathbb{Z}_{n_1}, V\subset \mathbb{Z}_{n_2}$ are intervals of lengths $1\leq |U|, |V|\leq b$. Then $J_0\cap V\not= \emptyset.$
\end{lem}

{\bf Proof.} Since $d(I_1, I_2)\geq b+1$, the projection $U$ can meet only one of them. Suppose $U\cap I_2=\emptyset.$Then the other projections of $R_2$ and $U\times V$ must meet: $J_0\cap V\not= \emptyset$ , as stated.
\qed

\smallskip

\begin{lem} Suppose that the positive integers $k, \ell, b, n_1, n_2$ satisfy the inequalities $k, \ell \leq b, 2(k+b)< n_1, 2(\ell+b)< n_2$.
Let $\R$ be a proj-intersecting family of $k\times \ell$  rectangles in $\mathbb{Z}_{n_1}\times \mathbb{Z}_{n_2}$. Suppose that $\R$ contains $\ell $ pieces of $b$-blocking pairs of form (2)
with distinct bases. Then there is a $\beta \in \mathbb{Z}_{n_2}$ such that the projection of every member of $\R$
contains it.
\end{lem}

{\bf Proof.}  Lemma 4 implies that the bases must pairwise intersect. Let the bases be $B_0, B_1, \ldots , B_{\ell -1}$ . By Lemma 1 $B_0, B_1, \ldots , B_{\ell -1}$ are $\ell$ consecutive intervals: all intervals of length $\ell$ in an interval $A$ of length $2\ell -1$. All these bases contain the middle element of $A$, this will play the role of $\beta$.

Apply Lemma 4 for an arbitrary member of $\R$ (that is $|U|=k, |V|=\ell.$) By this lemma $V$ intersects each $B_i$.
It is easy to see that $V$ must be equal to one of them. Hence it contains $\beta$. \qed

\smallskip
\begin{cor} $|\R|\leq \ell n_1$ holds under the conditions of Lemma 5.
\end{cor}

{\bf Proof.} $\ell n_1$ is the total number of $k\times \ell$ rectangles whose projection on $\mathbb{Z}_{n_2}$
contains a fixed element $\beta$.
\qed

\smallskip
\begin{lem} Suppose that the positive integers $k, \ell, b, n_1, n_2$ satisfy the inequalities $k, \ell \leq b, 2(k+b)< n_1, 2(\ell+b)< n_2$.
Let $\R$ be a proj-intersecting family of $k\times \ell$  rectangles in $\mathbb{Z}_{n_1}\times \mathbb{Z}_{n_2}$. Suppose that the number of $b$-blocking pairs of form (2)
with distinct bases is at least one and at most $\ell -1$. Then
$$|\R|\leq 4b^2+(\ell -1)n_1.$$
\end{lem}

{\bf Proof.} Define a partition of $\R^2=\R_{\rm b}^2\cup \R_{\rm s}^2$ where  $\R_{\rm b}^2$ denotes the set of members of $\R^2$ with multiplicity at least $k+b+1$. If $B\in \R_{\rm b}^2$ then, by Lemma 2 there is a $b$-blocking pair with basis $B$. One of the conditions of the lemma is that $|\R_{\rm b}^2|\leq \ell -1$. These members of $\R_{\rm b}^2$ have only the trivial bound on the multiplicities: $n_1$. We have
$$\sum_{B\in \R_{\rm b}^2}\mu (B)\leq (\ell -1)n_1.\eqno(5)$$

Choose one member $B_0\in \R_{\rm b}^2$. Lemma 4 implies that $B_0$ must meet every member of $\R_{\rm s}^2$.
There are at most $2\ell -2$ intervals of length $\ell$ meeting $B_0$. This results in $|\R_{\rm s}^2|\leq 2\ell -2.$
Suppose $J\in \R_{\rm s}^2$. By Lemma 2, if $\mu (J) \geq k+b+1$ then there is a $b$-blocking pair with basis $J$.
Therefore $\mu (J) \leq k+b$ holds. We obtained
$$\sum_{J\in \R_{\rm s}^2}\mu (J)\leq (2\ell -2)(k+b)\leq (2b)(2b)=4b^2.\eqno(6)$$
We can finish the proof using (5) and (6):
$$|\R|=\sum_{J\in \R^2}\mu (J) =\sum_{B\in \R_{\rm b}^2}\mu (B) + \sum_{J\in \R_{\rm s}^2}\mu (J)\leq 4b^2+(\ell -1)n_1.$$
\qed

The statements of Lemmas 3, 5 and 6 can be summarized in the following Corollary.

\begin{cor} Suppose that the positive integers $k, \ell, b, n_1, n_2$ satisfy the inequalities $k, \ell \leq b, 2(k+b)< n_1, 2(\ell+b)< n_2$.
Let $\R$ be a proj-intersecting family of $k\times \ell$  rectangles  in $\mathbb{Z}_{n_1}\times \mathbb{Z}_{n_2}$. Then one of the followings hold.
$$|\R|< 9b^2,$$
$$|\R|\leq 4b^2+(\ell -1)n_1,$$
$$|\R|\leq \ell n_1,$$
$$|\R|\leq 4b^2+(k-1)n_2,$$
$$|\R|\leq k n_2.$$
\end{cor}

\smallskip
This corollary, however is not sufficient for our final goal when there are rectangles of different sizes. The reason is that the statements of Corollary 2 cannot be independently used for different sizes, because they strongly interact. See the lemma below.

\begin{lem} Suppose that the positive integers $k_1, k_2, \ell_1, \ell_2, b, n_1, n_2$ satisfy the inequalities $k_1, k_2, \ell_1, \ell_2 \leq b, 4b< n_1, 4b< n_2$.
Let $\R_i (i=1,2)$ be a family of $k_i\times \ell_i$  rectangles  in $\mathbb{Z}_{n_1}\times \mathbb{Z}_{n_2}$ and suppose that $\R=\R_1\cup \R_2$ is proj-intersecting. Then $\R$ cannot simultaneously contain a $b$-blocking pair with basis in  $\mathbb{Z}_{n_1}$ (of the form (2)) and another one with basis in $\mathbb{Z}_{n_2}$ (of the form (3)).
\end{lem}

{\bf Proof.} Let $R_1, R_2\in \R_1$ be a $b$-blocking pair with basis $J_0\subset \mathbb{Z}_{n_2}$. In other words
$R_1= I_1\times J_0, R_2=I_2\times J_0$ where $d(I_1,I_2)\geq b+1$. On the other hand let $R_3, R_4\in \R_2$ be the other $b$-blocking pair with basis $I_0\subset \mathbb{Z}_{n_2}$ where
$R_3= I_0\times J_3, R_4=I_0\times J_4$ with $d(J_3,J_4)\geq b+1$.

Here $d(I_1, I_2)\geq b+1$ implies that $I_0$ can intersect only one of them. Suppose
$$I_0\cap I_2=\emptyset . \eqno(7)$$
On the other hand, using $d(J_3,J_4)\geq b+1$, the interval $J_0$ can meet at most one of $J_3$ and $J_4$. Suppose
$$J_0\cap I_4=\emptyset.\eqno(8)$$
(7) and (8) show that $R_2$ and $R_4$ are not proj-intersecting. This contradiction finishes the proof. \qed

\begin{lem} Suppose that the positive integers $k_i, \ell_i, b, n_1, n_2$ satisfy the inequalities $k_i, \ell_i \leq b (1\leq i\leq m), 4b< n_1, n_2$.
Let $\R_i$ be a family of $k_i\times \ell_i$  rectangles  in $\mathbb{Z}_{n_1}\times \mathbb{Z}_{n_2} (1\leq i\leq m)$. Suppose that $\R =\bigcup_{i=1}^m \R_i$ is a proj-intersecting family. Assume that there is a $b$-blocking pair say in $R_1$ with basis in $\mathbb{Z}_{n_2}$. Then either
$$|\R_i|< 9b^2\eqno(9)$$
or
$$|\R_i|\leq 4b^2+(\ell_i -1)n_1\eqno(10)$$
or
$$|\R_i|\leq \ell_i n_1\eqno(11)$$
holds for all $i (1\leq i\leq m).$
\end{lem}

{\bf Proof.} If (9) does not hold for an $i$ then by Lemma 3 there is a $b$-blocking pair in $\R_i$. By Lemma 7
it must be one with a basis in $\mathbb{Z}_{n_2}$. Suppose that the number of distinct bases of $b$-blocking pairs
$\R_i$ is $\ell_i$. Then by Lemma 5 and Corollary 1 we obtain (11). On the other hand if the number of distinct bases is between 1 and $\ell_i -1$ then Lemma 6 implies (10). \qed

\begin{rem} Of course, if the bases of the $b$-blocking pairs are in $\mathbb{Z}_{n_1}$ then (10) and (11) are replaced by
$$|\R_i|\leq 4b^2+(k_i-1)n_2\eqno(12)$$
or
$$|\R_i|\leq k_i n_2.\eqno(13)$$
\end{rem}

\begin{rem} The bases of $b$-blocking pairs for distinct $\R_i$'s must also intersect. This fact gives
a stricter structure for the system of bases, but we will not use this fact.
\end{rem}

Now we get rid of the cases, assuming that $n_1$ and $n_2$ are large enough.

\begin{lem} Suppose that the positive integers $k_i, \ell_i, b, n_1, n_2$ satisfy the inequalities $k_i, \ell_i \leq b (1\leq i\leq m), 9b^2< n_1, n_2$.
Let $\R_i$ be a family of $k_i\times \ell_i$  rectangles  in $\mathbb{Z}_{n_1}\times \mathbb{Z}_{n_2} (1\leq i\leq m)$. Suppose that $\R =\bigcup_{i=1}^m \R_i$ is a proj-intersecting family.  Then either
$$|\R_i|\leq \ell_i n_1\eqno(14)$$

holds for all $i (1\leq i\leq m)$ or
$$|\R_i|\leq k_i n_2. \eqno(15)$$
\end{lem}

{\bf Proof.} We only have to notice that each of (9) and (10) implies (11) under the condition $9b^2<n_1, n_2$,
while (12) implies (13).
\qed

\begin{cor} Let $\lambda_i>0 (1\leq i\leq m)$ be real numbers.
Under the conditions of Lemma 9
$$\sum_{i=1}^m\lambda_i |\R_i|\leq \max \left\{ n_1\sum_{i=1}^m\lambda_i \ell_i, n_2\sum_{i=1}^m \lambda_i k_i \right\}$$
holds.
\end{cor}

{\bf Proof.} Indeed, summing up (14) for $i$ in Lemma 9
 $$\sum_{i=1}^m\lambda_i |\R_i|\leq  n_1\sum_{i=1}^m\lambda_i \ell_i$$
 is obtained while (15) leads to
 $$\sum_{i=1}^m\lambda_i |\R_i|\leq  n_2\sum_{i=1}^m \lambda_i k_i.$$
 Since one of them must hold, the statement of the lemma follows. \qed

\smallskip
{\bf Proof of Theorem 3} Define the families
$$\F_i=\{ F\in \F:\ |F\cap X_1|=k_i, |F\cap X_2|=\ell_i\}.$$
We use double counting for the sum
$$\sum_{F, \C_1 ,\C_2}s(F)\eqno(16)$$
where $\C_j$ is a cyclic permutation of $\mathbb{Z}_{n_j} (j=1,2)$, $F\in \F$ and it forms a
rectangle for the product of these two cyclic permutations and the weight $s(F)$ is defined in the following way:
$$s(F)=s_i(F)={1\over n_1!}\cdot {1\over n_2!}{n_1\choose k_i}{n_2\choose \ell_i} {\mbox{ if }} F\in \F_i.\eqno(17)$$

For a fixed set $F\in \F_i$ there are $k_i!(n_1-k_i)!\ell_i !(n_2-\ell_i )!$ pairs of cyclic permutations  $(\C_1,\C_2 )$ in which $F$ is a rectangle.  That is (16) is equal to
$$\sum_{F\in \F} \sum_{ \C_1 ,\C_2}s(F)=\sum_{i=1}^m \sum_{F\in \F_i} \sum_{\C_1 ,\C_2}s_i(F)=$$   $$\sum_{i=1}^m|\F_i|k_i!(n_1-k_i)!\ell_i !(n_2-\ell_i )!{1\over n_1!}\cdot {1\over n_2!}{n_1\choose k_i}{n_2\choose \ell_i}=\sum_{i=1}^m|\F_i|=|\F|. \eqno(18)$$
Now fix the permutations in (16):
$$\sum_{\C_1 ,\C_2}\sum_Fs(F)=\sum_{\C_1 ,\C_2}\sum_{i=1}^m|\R_i|s_i(F)\eqno(19)$$
where $\R_i$ denotes the set of rectangles obtained from $\F_i$ in these fixed cyclic permutations. $\F$ is an intersecting family. It is easy to see that this implies that $\R=\cup \R_i$ is proj-intersecting.  Corollary 3
can be applied with $\lambda_i=s_i(F)$ since its conditions are satisfied.

$$\sum_{i=1}^m|\R_i|{1\over n_1!}\cdot {1\over n_2!}{n_1\choose k_i}{n_2\choose \ell_i} \leq $$
$$\max \left\{ n_1\sum_{i=1}^m{1\over n_1!}\cdot {1\over n_2!}{n_1\choose k_i}{n_2\choose \ell_i}\ell_i,
n_2\sum_{i=1}^m{1\over n_1!}\cdot {1\over n_2!}{n_1\choose k_i}{n_2\choose \ell_i}k_i \right\} =$$
$$\max \left\{ \sum_{i=1}^m{1\over (n_1-1)!}\cdot {1\over (n_2-1)!}{n_1\choose k_i}{n_2-1\choose \ell_i-1},\right. \hskip 50mm$$
$$\hskip 48mm\left. \sum_{i=1}^m{1\over (n_1-1)!}\cdot {1\over (n_2-1)!}{n_1-1\choose k_i-1}{n_2\choose \ell_i} \right\} .
 \eqno(20)$$
 Taking into account that the number of pairs of cyclic permutations $\C_1, \C_2$ is $(n_1-1)!(n_2-1)!$ and using (20) we obtain an upper estimate on (19):
 $$\max \left\{ \sum_{i=1}^m{n_1\choose k_i}{n_2-1\choose \ell_i-1},
\sum_{i=1}^m{n_1-1\choose k_i-1}{n_2\choose \ell_i} \right\} .
 \eqno(21)$$
 Since both (18) and (19) are equal to (16), we arrived to the inequality
 $$|\F|\leq \max \left\{ \sum_{i=1}^m{n_1\choose k_i}{n_2-1\choose \ell_i-1},
\sum_{i=1}^m{n_1-1\choose k_i-1}{n_2\choose \ell_i} \right\}.$$
The quantities in the $\max$ are the numbers of all sets containing a fixed element of $X_2$ and $X_1$, respectively. \qed

\section{Open problems}
The upper bound in the Erd\H os-Ko-Rado theorem is reached for the family of all $k$-element sets containing a fixed element. What happens if we exclude this construction?
Hilton and Milner found the largest intersecting, but not trivially intersecting family.

\begin{thm} {\rm{\cite{HM}}}
If $\F$ is an intersecting but not a trivially intersecting family, $\F \subset {[n]\choose k} (2k\leq n)$ then
$$|\F|\leq 1+{n-1\choose k-1}-{n-k-1\choose k-1}.$$
\end{thm}

The construction giving equality is the following. Fix an element $x\in X$ and a $k$-element set $K$ such that $x\not\in K, K\subset X$. The extremal family will consist of all $k$-element sets containing $x$ and intersecting $K$.

This construction can be imitated for two parts, but the size of the family depends on weather $x$ and $K$ are in $X_1$ or $X_2$. Suppose that $x\in X_1, K\subset X_1$. Then the construction of the family is the following:
$$\F =\{ F:\ x\in F, F\cap K\not=\emptyset , |F\cap X_1|=k, |F\cap X_2|=\ell \} .$$ The case when $x$ and $K$ are in $X_2$ is analogous.
We conjecture that one of these constructions is the best.

\begin{conj}
 If $\F$ is a non-trivially intersecting subfamily ${X_1, X_2\choose k,\ell}$
then
 $$|\F|\leq
\max \left\{ \left( 1+{n_1-1\choose k-1}-{n_1-k-1\choose k-1}\right){n_2\choose \ell},\right. \quad \quad \quad \quad \quad \quad \quad \quad $$
  $$ \quad \quad \quad \quad \qquad \qquad \quad \quad  \qquad \qquad \left.{n_1\choose k}\left( 1+{n_2-1\choose \ell-1}-{n_2-\ell-1\choose \ell-1}\right) \right\} . $$
\end{conj}
The constructions in this conjecture are non-trivially intersecting, but the ``intersections happen" in one side.
Our next question is
what happens if one side is not enough for satisfying the intersection conditions. The families
satisfying the following conditions are called {\it two-sided intersecting}:
there are members $F_{11}, F_{12}, F_{21}, F_{22}\in \F$ such that
$F_{11}\cap F_{12}\cap X_1=\emptyset$ and $F_{21}\cap F_{22}\cap X_2=\emptyset$ hold.

To better understand our ``best" two-sided intersecting construction one more notion and one more statement are needed.
The families $\F, \G  \subset {[n]\choose k}$ are {\it cross-intersecting}
if $F\cap G$ is non-empty for every pair of members $F\in \F, G\in \G$. Here the total number $|\F|+|\G|$ of members should be maximized. But only for non-empty families, otherwise it is trivial and uninteresting.

\begin{thm} {\rm{\cite{HM}}}
If $\F, \G  \subset {[n]\choose k}$ are non-empty
cross-intersecting families then
$|\F|+|\G|\leq 1+{n\choose k}-{n-k\choose k}$.
\end{thm}

This estimate is sharp: let $K$ be a $k$-element subset, let $\F_1(n) =\{ K\}$ and let $\G_1(n)$ consist of all $k$-element
subsets intersecting $K$.

Now we are ready to construct a large two-sided intersecting family $\F$.
It will be done in two steps. First the ``projection" $\F_2$ of our family $\F$ for $X_2$ will be given,
formally
$$\F_2=\{ F\cap X_2:\ F\in \F \}.\eqno(22)$$
Then we will determine the families $\F(M)$ ``belonging" to the members $M\in \F_2$:
$$\F(M)=\{ F\subset X_1:\ F\cup M\in \F \}.$$
$\F_2$ will be an ``almost intersecting" family in which there is only one non-intersecting pair.
Start with a family extremal for Theorem 4 in $X_2$ where $x$ is a fixed element, $L$ is a
fixed $\ell$-element subset, not containing $x$:
$$\{ F:\ x\in F, F\cap L\not=\emptyset\}.$$
Add another $\ell$-element set $L^{\prime}$ such that $x\in L^{\prime}$ holds.
Of course, $L\cap L^{\prime}=\emptyset$.
Define the ``projection" (22) as
$$\F_2= \{ F:\ x\in F, F\cap L\not=\emptyset\}\cup \{ {L^{\prime}} \}.$$
$L\cap L^{\prime}=\emptyset,$ all other pairs are intersecting.
Hence $\F(M)$ can be chosen to be ${X_1\choose k}$ if $M\not=L$ or $L^{\prime}$. However the families
$\F(L)$ and $\F(L^{\prime})$ must be a pair of non-empty cross-intersecting families. To maximize
the sum of their sizes the construction of Theorem 5 should be used. Choose a $k$-element subset
$K\subset X_1$ and define $\F(L)=\{ K\}, \F(L^{\prime})=\{ F\subset X_1:\ F\cap K\not= \emptyset \}.$
We believe that either this, or its symmetric version is the largest such family.

\begin{conj}
 If $\F$ is a two-sided intersecting subfamily of ${X_1, X_2\choose k,\ell}$
then
 $$|\F|\leq
\max \left\{ \left( {n_2-1\choose \ell -1}-{n_2-\ell-1\choose \ell-1}\right){n_1\choose k}+1+{n_1\choose k}-{n_1-k\choose k},\right. \quad \quad \quad \quad \quad \quad \quad \quad $$
  $$ \quad \quad \quad \qquad \left.  \left( {n_1-1\choose k -1}-{n_1-k-1\choose k-1}\right){n_2\choose \ell}+1+{n_2\choose \ell}-{n_2-\ell \choose \ell} \right\} . $$
\end{conj}

\section{Related results from the past}

{\bf
$k$-part Sperner theorems.} Let $\F$ be a family of subsets of the $n$-element $X$ satisfying the condition
 that $F\in \F, G\in \F$ implies $F\not\subset G$. (In other words $\F$ is inclusion-free.) The well-known theorem
of Sperner \cite{Sp} states that $|\F|\leq {n\choose \lfloor n/2 \rfloor}$ and this estimate is sharp: take all sets of size $\lfloor n/2 \rfloor$.

Let now $X$ be a disjoint union of $X_1$ and $X_2$ where $|X_1|+|X_2|=n.$ It was noticed in \cite{K2Sp} and \cite{Kl}
that Sperner's bound remains valid if only such pairs of subsets $F\subset G$ are excluded which are equal in one of the parts, that is either $F\cap X_1=G\cap X_1$ or $F\cap X_2=G\cap X_2$ holds.

But this statement does not generalize to three parts. One can give a family $\F$ containing more than ${n\choose \lfloor n/2 \rfloor}$ subsets of $X=X_1\cup X_2 \cup X_3$ such that there are no pairs $F\in \F, G\in \F, F\subset G$ with
$G-F\subset X_i$ for some $i$. It seems to be difficult to determine the largest $\F$ under this condition. The best
estimates can be found in \cite{M}.

{\bf Non-empty intersections in both parts.} Sali \cite{Sa} considered the following problem. Let
$$\F \subset {X_1, X_2\choose k,\ell}$$
be a $t+1$-intersecting family with an additional condition: $F\cap G\cap
X_1\not= \emptyset, F\cap G\cap X_2\not= \emptyset$
for all $F,G \in \F$. The exact maximum
size of $|\F|$ is determined under this condition if both $|X_1|$ and
$|X_2|$ are large enough. In spite of formal similarities this problem
is very different from our problem in nature.

{\bf Using the two parts differently.} The following result is also slightly related. Let
$$\F \in {X_1\choose k}\cup {X_1\cup X_2\choose \ell}$$
that is the members of the family either have exactly $k$ elements in the first part, or exactly
$\ell$ elements in the whole underlying set. Wang and Zhang \cite{WZ} determined the maximum size of an intersecting family of this form.

\section{Acknowledgement} The author is indebted to Noga Alon, Peter Frankl and Casey Tompkins for helpful discussions and the anonymous referee for careful reading and corrections.

\end{document}